\documentclass[leqno,12pt]{amsart} 
\usepackage{amssymb} 
\theoremstyle{plain} %text of this environment is typesetted in italics
\newtheorem{theorem}{\indent\sc Theorem}[section]
\newtheorem{lemma}[theorem]{\indent\sc Lemma}
\newtheorem{corollary}[theorem]{\indent\sc Corollary}
\newtheorem{proposition}[theorem]{\indent\sc Proposition}

\theoremstyle{definition} %text of this environment is typesetted in roman letters

\newtheorem{remark}[theorem]{\indent\sc Remark}

\numberwithin{equation}{section}

\def\<{\left < }
\def\>{\right >}
\def\({\left ( }
\def\){\right )}

   \def\p{\partial}

\def\sn{{\rm sn}}
\def\cn{{\rm cn}}
\def\dn{{\rm dn}}

\begin{document}

\title[Lagrangian surfaces in complex Euclidean plane] {Lagrangian surfaces in complex Euclidean plane via spherical and hyperbolic curves}
\author[I. Castro ]{Ildefonso Castro}

\thanks{First author's research is partially supported  by a MEC--FEDER grant No. MTM2004-00109.}

\author[B.-Y. Chen]{Bang-Yen Chen}
\thanks{A part of this article was done while the second author
visited Universidad de Ja\'en, Spain. The second author would like to
express his many thanks for the hospitality he received during his
visit.}

\begin{abstract}

We present a method to construct a large family of
Lagrangian surfaces in complex Euclidean plane ${\bf C}^2$ by using
Legendre curves in the 3-sphere and in the anti de Sitter 3-space
or, equivalently, by using spherical and hyperbolic curves, respectively.
Among this family, we characterize minimal, constant mean curvature, Hamiltonian-minimal and
Willmore surfaces in terms of simple properties of the curvature of
the generating curves. As applications, we provide explicitly conformal parametrizations of known and new
examples of  minimal, constant mean curvature, Hamiltonian-minimal and Willmore surfaces in
${\bf C}^2$.

\end{abstract}



\thanks{}

\maketitle

\section{Introduction}

An immersion $\phi : M^n \to \tilde M^n$ of an $n$-manifold $M^n$
into a Kaehler $n$-manifold $\tilde M^n$ is called a {\em Lagrangian
immersion} if the complex structure $J$ of $\tilde M^n$
interchanges each tangent space of $M^n$ with its corresponding normal space.
Lagrangian submanifolds appear naturally in several contexts of
mathematical physics. A very important problem in this setting is to
find nontrivial examples of Lagrangian submanifolds with some given
geometric properties.
In this line, we find many papers (see the survey article \cite{c})
where the different authors investigate intrinsic and extrinsic
geometric properties related mainly with the intrinsic curvatures and
the mean curvature vector of the submanifolds, respectively.

An important problem in the theory of Lagrangian
surfaces is to find non-trivial examples  with some given
special  geometric properties.
For instance,  a method
was given in \cite{J1}  to construct an important
family of special Lagrangian submanifolds in {\bf C}$^n$ with large symmetric groups.  Also,  a spinor-like representation formula for Lagrangian surfaces in ${\bf C}^2$ which  parameterizes  immersions through two complex functions $F_1,F_2$ and a real one (the Lagrangian angle) was introduced in  \cite{a1}. This  formula is  useful to construct  examples of  Lagrangian surfaces in {\bf C}$^2$.

In this article, we present  a simple specific new method to construct
Lagrangian surfaces in complex Euclidean plane ${\bf C}^2$ with nice
properties  that only involves two Legendre curves; one in the
3-sphere $S^3$ and the other in the anti De Sitter 3-space $H^3_1$.

Recall that a regular curve $\gamma: I_1 \to S^3 $
(respectively $\alpha: I_2 \to H^3_1 $) is called a {\em Legendre
curve} if $\< \gamma ' (s), i \gamma (s) \> =0$ (respectively $\<
\alpha ' (t), i \alpha (t) \> =0$) holds identically. For each such
pair of Legendre curves $(\gamma,\alpha)$, we consider the map:
\begin{align}\label{1.1}\phi:I_1\times I_2\to {\bf C}^2
={\bf C}\times{\bf C};\;\; (t,s)\mapsto
\(\alpha_1(t)\gamma_1(s),\alpha_2(t)\gamma_2(s)\).
\end{align}

In Section 2, we show that the map $\phi$ defined by \eqref{1.1} is a
Lagrangian immersion. We also study geometric properties of such Lagrangian
surfaces. In particular, we investigate the close relationship of
such a Lagrangian surface with the curve in the 2-sphere $S^2$ and the
curve in the hyperbolic 2-plane $H^2$ given by the projections of $\gamma $ and
$\alpha $ via their corresponding Hopf fibrations. 

In Section 3 we  prove an useful ``additive formula'' (see Theorem 1) involving the
Lagrangian angle map of $\phi $  and  the Legendre angles of $\gamma$ and $\alpha$ for the Lagrangian immersions. As a consequence, we establish a simple 
relationship (see Corollary 1) between the mean
curvature vector of $\phi $ and the curvature functions of $\gamma $
and $\alpha$. 

 The last section provides several nice applications of the results obtained in Section 3. 
First we characterize the minimal
Lagrangian surfaces obtained by our construction in terms of
geodesics in $S^3$ and $H^3_1$. In such a way we are able to
provide explicit expressions  of the minimal Lagrangian conformal
immersions  in ${\bf C}^2$ in terms of some elementary functions.
We also determine Lagrangian surfaces constructed by our method with constant mean curvature and, in particular,  with parallel mean curvature vector. This enables us to obtain interesting new examples of Lagrangian tori in ${\bf C}^2$ with constant mean curvature. Next we
 characterize  Hamiltonian-minimal Lagrangian
surfaces among the family of Lagrangian surfaces using our
construction with Legendre curves such that their curvature functions (in terms of the arclength parameter)
are linear. 
As a by-product, we are able to establish the  explicit expressions of some new
Hamiltonian-minimal Lagrangian conformal immersions in ${\bf C}^2$ in
terms of elementary functions as well.
Finally,  we apply our construction method to provide new
examples of Willmore Lagrangian surfaces in ${\bf C}^2$. Our result
states that the Lagrangian surfaces constructed by the pair
$(\gamma,\alpha)$ of Legendre curves are Willmore surfaces
if and only if  $\gamma $ and $\alpha $ are elastic curves in $S^3$
and $H^3_1$, respectively.

\section{A new construction method of Lagrangian surfaces}

In the complex Euclidean plane ${\bf C}^2$ we consider the bilinear
Hermitian product defined by
\[
(z,w)=z_1\bar{w}_1+z_2\bar{w}_2, \quad  z,w\in{\bf C}^2.
\]
Then $\langle\, \, , \, \rangle = \Re (\,\, , \,)$ is the Euclidean
metric on ${\bf C}^2$ and
$\omega = -\Im (\,,)$ is the Kaehler two-form given by $\omega
(\,\cdot\, ,\,\cdot\,)=\langle J\cdot,\cdot\rangle$, where $J$ is the
complex structure on ${\bf C}^2$. 

Let $\phi:M\to {\bf C}^2$ be an isometric Lagrangian immersion of a
surface $M$ into ${\bf C}^2$, i.e. an immersion satisfying
$\omega_{|M}\equiv 0$.
We denote the Riemannian connections of $M$ and ${\bf C}^2$ by
$\nabla$ and $\bar{\nabla}$, respectively. We also denote by $\langle
\,\, , \, \rangle $ the induced metric on $M$.
Then we have $\phi^* T{\bf C}^2 =\phi_* TM \oplus T^\perp M$, where
$TM$ and $T^\perp M$ are the tangent and normal bundles of $M$,
respectively.
The second fundamental form $\sigma $ is given by $\sigma
(x,y)=JA_{Jx}y$, where $A$ is the shape operator. Thus, the
(0,3)-tensor $C(x,y,z)=\left<\sigma(x,y),Jz\right>$ is totally
symmetric.

The space of oriented Lagrangian planes in ${\bf C}^2$ can be
identified with the symmetric space $U(2)/SO(2)$, so
the determinant map, $\det:U(2)/SO(2)\to S^1$, is well-defined. If
$M$ is an orientable Lagrangian
surface in ${\bf C}^2$ and $\nu:M\to U(2)/SO(2)$ is its Gauss map,
then $\det\circ\nu:M\to S^1$ can be expressed as
$\det\circ\nu=e^{i\beta_\phi}$ for some function $\beta_\phi:M\to
{\bf R} /2\pi{\bf Z}$. This function $\beta_\phi$ is called the {\em
Lagrangian angle map} of $\phi$.
The Lagrangian angle map $\beta_\phi$ satisfies
\begin{align}\label{beta}
J\nabla\beta_\phi=2H,
\end{align}
where $H$ is the mean curvature of $\phi$,
defined by $H=(1/2) {\rm trace} \, \sigma$.

The Lagrangian immersion $\phi $ is called minimal if $H=0$
identically, or equivalently, the Lagrangian angle map $\beta_\phi$
is constant. The minimality on the surfaces means that the surface is a
critical point  of the area functional over any compactly supported variation. 
On
the other hand, Hamiltonian-minimal Lagrangian surfaces are
Lagrangian surfaces which are critical points of the area functional
with respect to a special class of infinitesimal variations
preserving the Lagrangian constraint;  namely, the class of compactly
supported Hamiltonian vector fields  (see \cite{O}). Such Lagrangian
surfaces are characterized by the harmonicity of the Lagrangian angle
map $\beta_\phi$ (cf. \cite{HR}).

\vspace{0.2cm}

Let $S^3$ and $H^3_1$ denote the
unit hypersphere  and the unit anti De Sitter space in
${\bf C}^2$ given respectively by
$$S^3=\left\{ (z,w)\in {\bf C}^2 ; |z|^2+|w|^2=1 \right\}
,\;\;H^3_1=\left\{ (z,w)\in {\bf C}^2 ; |z|^2-|w|^2=-1
   \right\} .$$

Let $\gamma:= \gamma(s)=(\gamma_1,\gamma_2): I_1 \to  S^3$ be a unit
speed Legendre curve in  $S^3$  and
$\alpha=\alpha(t)=(\alpha_1,\alpha_2):I_2\to H^3_1$  a unit speed
Legendre curve in $H^3_1$.
Then $\gamma $ and $\alpha $ satisfy
\begin{align}&\label{2.1} |\gamma_1|^2+|\gamma_2|^2=1,\hskip.2in
|\gamma'_1|^2+|\gamma'_2|^2=1,\;\;\, \gamma'_1
\bar\gamma_1+\gamma'_2\bar\gamma_2=0,\;\;
    \\&\label{2.2}|\alpha_1|^2-|\alpha_2|^2=-1,\;\;
|\alpha'_1|^2-|\alpha'_2|^2=1,\;\; \alpha'_1
\bar\alpha_1-\alpha'_2\bar\alpha_2=0.\end{align}

\begin{proposition} \label{P:1}
   Let $\gamma$ be a unit speed
Legendre curve in
$S^3$ and $\alpha$ be a unit speed
Legendre curve in $H^3_1$.
Consider the map: $\phi:I_1\times I_2\subset {\bf R}^2\to {\bf C}^2
={\bf C}\times{\bf C}$
defined by
\begin{align}\label{2.3}
\phi(t,s)=\(\alpha_1(t)\gamma_1(s),\alpha_2(t)\gamma_2(s)\).
\end{align}

Then $\phi$ is a Lagrangian conformal immersion in ${\bf C}^2$
such that the induced metric is given by
\begin{align}\label{2.4}\left<\,\;,\;\right>=\(|\gamma_1|^2
+|\alpha_1|^2\)\(dt^2+ds^2\)
\end{align}
and the intrinsic tensor $C(x,y,z)=\left<\sigma(x,y),Jz\right>$
is given by
\begin{equation}
\begin{aligned}\label{2.5}&C_{ttt}
   =\left<\alpha_1'',J\alpha_1'\right>|\gamma_1|^2
+\left<\alpha_2'',J\alpha_2'\right>|\gamma_2|^2,
\\&C_{tts}=\left<{\gamma}'_1,J\gamma_1\right>,
\\& C_{tss}=\left<\alpha_1',J\alpha_1\right>,\\&
C_{sss}=|\alpha_1|^2
\left<\gamma_1'',J\gamma_1'\right>+|\alpha_2|^2
\left<{\gamma}_2'',J\gamma_2'\right>,\end{aligned}\end{equation}
where  $C_{ttt}=C(\p_t,\p_t,\p_t),C_{tts}=C(\p_t,\p_t,\p_s),
\cdots$, etc.
The $J$ in \eqref{2.5} is the $+\pi/2$-rotation acting on ${\bf
C}\equiv {\bf R}^2$. \end{proposition}

\begin{proof} From \eqref{2.3} we get
\begin{align}\label{2.6}\phi_t=(\alpha'_1(t)\gamma_1(s),
\alpha'_2(t)\gamma_2(s)),\;\;
\phi_s=(\alpha_1(t)\gamma'_1(s),\alpha_2(t)\gamma'_2(s)).\;
\end{align}
Thus, by applying \eqref{2.1} and \eqref{2.2}, we find
\begin{equation}
\begin{aligned}\label{2.7} |\phi_t|^2&=|\alpha'_1|^2|\gamma_1|^2
+ |\alpha'_2|^2|\gamma_2|^2 \\&=|\alpha'_1|^2|\gamma_1|^2+
(|\alpha'_1|^2-1)(1-|\gamma_1|^2)
\\&=|\alpha'_1|^2+|\gamma_1|^2-1
\\&=|\alpha'_2|^2+|\gamma_1|^2 .
\end{aligned}\end{equation}

On the other hand, from the last equation of \eqref{2.2}, we have
\begin{equation}
\begin{aligned}\label{2.8} & |\alpha_1|^2(1+
|\alpha'_2|^2)=|\alpha_1|^2 |\alpha'_1|^2
\\&\hskip.2in  =|\alpha_2|^2 |\alpha'_2|^2=|\alpha'_2|^2(1+ |\alpha_1|^2).
\end{aligned}\end{equation}
Thus, we obtain $|\alpha_1|^2=|\alpha'_2|^2$. Substituting this into
\eqref{2.7} gives $|\phi_t|^2=|\alpha_1|^2+|\gamma_1|^2 $.
Similarly, we also have $|\phi_s|^2=|\alpha_1|^2+|\gamma_1|^2 $.
By the last equations in
\eqref{2.1} and \eqref{2.2}, we have
\begin{equation}
\begin{aligned}\label{2.9} & ( \phi_t,\phi_s )=\alpha'_1
\bar\alpha_1\gamma_1\bar\gamma'_1+\alpha'_2\bar\alpha_2
\gamma_2\bar\gamma'_2=0.
\end{aligned}\end{equation} Thus, by taking the imaginary part in
\eqref{2.9}, we see that $\phi $ is a Lagrangian immersion whose
induced metric via $\phi$ is given by \eqref{2.4}.

It follows from \eqref{2.1} and \eqref{2.3} that
\begin{equation}
\begin{aligned}\label{2.10}
( \phi_{tt}, \phi_s )&=( (\alpha''_1\gamma_1,
\alpha''_2\gamma_2),
(\alpha_1\gamma_1',\alpha_2\gamma'_2) ) \\&=
\alpha''_1
\bar\alpha_1\gamma_1\bar\gamma'_1 + \alpha''_2\bar\alpha_2
\gamma_2\bar\gamma'_2\\&=
\gamma_1\bar\gamma'_1(\alpha''_1
\bar\alpha_1-\alpha''_2\bar\alpha_2)\\&=-
\gamma_1\bar\gamma'_1,
\end{aligned}\end{equation}
where we have applied the identity: $\alpha''_1
\bar\alpha_1-\alpha''_2\bar\alpha_2=-1$ deduced from the last
equation of \eqref{2.2}. Similarly, we also have
\begin{equation}
\begin{aligned}\label{2.11}
(\phi_{tt},\phi_t)&=\alpha''_1\bar\alpha_1'|\gamma_1|^2+
\alpha''_2\bar\alpha'_2|\gamma_2|^2.
\end{aligned}\end{equation}
By taking the imaginary parts in \eqref{2.10} and \eqref{2.11}, we
obtain the first two equations of \eqref{2.5}. Similarly,
we also have the last two equations of \eqref{2.5}.
\end{proof}

In the same spirit as the proof  of \eqref{2.5},  we find by taking the real parts that
\begin{equation}\begin{aligned}\label{2.12}&\<\nabla_{\partial_t}\partial_t,
\partial_t\>=
\left<\alpha_1'',\alpha_1'\right>|\gamma_1|^2
+\left<\alpha_2'',\alpha_2'\right>|\gamma_2|^2,
\\&\<\nabla_{\partial_t}\partial_t,
\partial_s\>=-\<\nabla_{\partial_t}\partial_s,
\partial_t\>=-\left<\gamma_1',\gamma_1\right>,
\\&
\<\nabla_{\partial_t}\partial_s,\partial_s\>=-
\<\nabla_{\partial_s}\partial_s,\partial_t\>=\left<\alpha_1',
\alpha_1\right>,\\&
\<\nabla_{\partial_s}\partial_s,
\partial_s\>=|\alpha_1|^2
\left<\gamma_1'',\gamma_1'\right>+|\alpha_2|^2
\left<\gamma_2'',\gamma_2'\right>.\end{aligned}\end{equation}
Since $\phi$ is a conformal Lagrangian immersion, the
Laplacian of the Lagrangian surface with the induced metric
\eqref{2.4} is given by
\begin{align}\label{2.13} \Delta=e^{-2u}\(\frac{\partial^2}{\partial t^2}
+\frac{\partial^2}{\partial s^2}\),\end{align} where
$e^{2u}=|\gamma_1|^2+|\alpha_1|^2$.
\vspace{0.2cm}

Via the Hopf fibration, Legendre curves in $S^3$ and $H^3_1$ are
projected into  curves in $S^2$ and $H^2$, respectively. Hence,  it
is interesting to describe the geometry of the Lagrangian surfaces
obtained in Proposition \ref{P:1} by using the geometry of the
curves  in $S^2$ and $H^2$. We study this as follows:

Let  $S^2(1/2):=\{ (x_1,x_2,x_3)\in {\bf R}^3 ; x_1^2+x_2^2+x_3^2=1/4
\}$ which is the 2-sphere with radius $1/2$ in {\bf R}$^3$. The Hopf fibration $\pi:S^3\to
S^2(1/2)\equiv CP^1(4)$ is given by
\begin{align}\label{3.1}
\pi(z,w)= \frac{1}{2} \(2z\bar w, |z|^2-|w|^2 \),\;\; (z,w)\in S^3\subset {\bf C}^2.
\end{align}
Notice that \eqref{3.1} is well-defined, since $|2z\bar
w|^2+(|z|^2-|w|^2)^2=(|z|^2+|w|^2)^2=1$.

For each Legendre curve $\gamma=\gamma(s)$ in $S^3$, the
projection $\xi=\pi\circ \gamma$ is a curve in $S^2(1/2)$.
Conversely,   each curve $\xi$ in $S^2(1/2)$ gives rise to a
horizontal lift $\tilde{\xi}$ in $S^3$ via $\pi$ which is unique up to a factor
$e^{i\theta_1},\theta_1\in {\bf R}$. Notice that each horizontal lift
of $\xi$  is a Legendre curve in $S^3$.

Since the Hopf fibration $\pi $ is a Riemannian submersion, each
unit speed Legendre curve $\gamma $ in $S^3$  is projected to   a unit speed curve $\xi$ in  $S^2(1/2)$ with the same curvature function.  From \eqref{3.1}, it is not difficult to see that
\begin{align}\label{flasCP}
|\gamma_1|^2=\frac{1}{2}+\xi_3, \, \,
\langle \gamma_1',J\gamma_1\rangle=(\xi \times \xi')_3,
\end{align}
where $\times $ denotes the cross product in ${\bf R}^3$
and $(\xi \times \xi')_3$ is the third coordinate of $\xi \times
\xi'$ in the 3-space ${\bf R}^3$ containing $S^2(1/2)$.

Similarly, let  $H^2(-1/2)=\{ (x_1,x_2,x_3)\in {\bf R}^3 ;
x_1^2+x_2^2-x_3^2=-1/4, \, x_3 \geq 1/2  \}$ which is the model of
the real hyperbolic plane of curvature $-4$. The Hopf fibration
$\pi:H^3_1\to H^2(-1/2)\equiv CH^1(-4)$ is then given by
\begin{align}\label{3.1bis}
\pi(z,w)= \frac{1}{2} \(2z\bar w, |z|^2+|w|^2 \),
\;\; (z,w)\in H^3_1\subset {\bf C}^2_1. \end{align}
Notice that \eqref{3.1bis} is well-defined, since $|2z\bar
w|^2-(|z|^2+|w|^2)^2=-(|z|^2-|w|^2)^2=-1$.

For each Legendre curve $\alpha=\alpha(t)$ in $H^3_1$, the
projection $\eta=\pi\circ \alpha$ is a curve in $H^2(-1/2)$.
Conversely,  each curve $\eta$
in $H^2(-1/2)$ gives rise to a horizontal lift $\tilde{\eta}$ in
$H^3_1$ via $\pi $ that is unique up to a factor
$e^{i\theta_2},\theta_2\in {\bf R}$. Each horizontal lift
$\tilde{\eta}$ is a Legendre curve in $H^3_1$.

In the same way as $\gamma$, if $\alpha $ is a unit speed Legendre
curve in $H^3_1$, then the projection
$\eta $ is also a unit speed curve in  $H^2(-1/2)$ with the same
curvature function. It follows from \eqref{3.1bis} that
\begin{align}\label{flasCH}
|\alpha_1|^2=-\frac{1}{2}+\eta_3, \, \,
\langle \alpha_1',J\alpha_1\rangle=(\eta \times \eta')_3.
\end{align}

Taking the above considerations into account, our construction of the
Lagrangian conformal surfaces in Proposition \ref{P:1} can also be
obtained by using a unit speed curve $\xi $ in $S^2(1/2)$
and a unit speed curve $\eta $ in $H^2(-1/2)$ as follows:
\begin{align}\label{2.3bis}
\phi(t,s)=\(\tilde\eta_1(t)\tilde\xi_1(s),\tilde\eta_2(t)\tilde\xi_2(s)\).
\end{align}

Notice that if we choose different horizontal lifts, say $\hat\xi $
and $\hat \eta $  of  $\xi $ and $\eta $, then we have $\hat\xi =
e^{i\theta_1}\tilde\xi$ and $\hat\eta = e^{i\theta_2}\tilde\eta$ for
some $\theta_1,\theta_2\in{\bf R}$. Hence, the corresponding
Lagrangian conformal immersion
$$
\psi(t,s)=\(\hat\eta_1(t)\hat\xi_1(s),\hat\eta_2(t)\hat\xi_2(s)\)
$$
is related with \eqref{2.3bis} by $\psi =
e^{i(\theta_1+\theta_2)}\phi $. Therefore,  the two Lagrangian
conformal immersions $\phi$ and $\psi$ via different horizontal lifts
are always congruent.

In fact, the geometry of the Lagrangian conformal immersion
$\phi $  depends essentially on the initial curves $\xi $ and $\eta
$. For example, it follows from \eqref{flasCP} and \eqref{flasCH}
that the induced metric of $\phi$ is given by
$ \left<\,\;,\;\right>=(\eta_3(t) + \xi_3(s))\(dt^2+ds^2\)$
and the intrinsic tensor $C$ satisfies
$C_{tts}=(\xi \times \xi')_3$, $C_{tss}=(\eta \times
\eta')_3,\ldots$, etc. These  show that the third coordinates of the
position of $\xi $ and $\eta$ in ${\bf R}^3$ are relevant in the
geometry of our construction. More precisely, any rotation around the
$x_3$-axis of ${\bf R}^3$ acting on the generating curves $\xi $ and
$\eta$ gives rise to a congruent Lagrangian immersion, since we have
\begin{align}\label{rot1}
{(e^{i\varphi_1}(\xi_1+i\xi_2);\xi_3)}=\pi(e^{i\varphi_1}\tilde\xi_1,\tilde\xi_2)
\end{align}
and
\begin{align}\label{rot2}
{(e^{i\varphi_2}(\eta_1+i\eta_2);\eta_3)}=\pi(e^{i\varphi_2}\tilde\eta_1,\tilde\eta_2).
\end{align}

As an illustration, the totally geodesic Lagrangian planes can be
obtained by taking any meridian in
$S^2(1/2)$ passing through the north and south poles and any meridian passing through
the vertex $(0,0,1/2)$ in
$H^2(-1/2)$. Therefore, up to congruence,  the totally geodesic Lagrangian
planes can be given  by $\phi(t,s)=(\cos s \sinh t, \sin s \cosh t)$.

\section{Additive formula of Lagrangian angle map}

For a unit speed Legendre curve
$\gamma=(\gamma_1,\gamma_2)$ in
$S^3$, we define the {\em Legendre angle}  $\theta_\gamma$ of
$\gamma$ by
\begin{align} \label{5.1}
e^{i\theta_\gamma}=\det{}_{\bf C}(\gamma,\gamma')=\gamma_1
\gamma'_2-\gamma'_1 \gamma_2 . %\;\; ({\rm mod}\; 2\pi).
\end{align}
For instance, the Legendre angle of
$\gamma(s)=(\cos s,\sin s)$ is 0 (mod $2\pi$).

\begin{lemma}  \label{L:1} Let $\gamma:I_1\to S^3\subset
{\bf C}^2$ be a unit speed curve. Then we have
\begin{enumerate}
\item[(1a)] If $\gamma$ is a Legendre curve, then it is a solution of the second order
differential equation$:$
\begin{align}\label{3.2} \gamma''-ik_\gamma \gamma'+\gamma=0, \end{align}
where $k_\gamma$ is the curvature of $\gamma$ in $S^3$.
\item[(1b)] If $\gamma $ satisfies \eqref{3.2}, then $\gamma $ is a Legendre curve if and only if 
$\left<\gamma' (0),i\gamma (0)\right>=0$ $(0\in I_1)$.
\item[(2)] The Legendre angle $\theta_\gamma$ satisfies $\theta'_\gamma=k_\gamma$.
\end{enumerate}\end{lemma}

\begin{proof} Statement (1a) has been proved in \cite{c1}. Statement (1b) follows from the constancy of the function
$s\mapsto \left<\gamma' (s),i\gamma (s)\right>$ using \eqref{3.2}.
Finally,  from \eqref{3.2} and \eqref{5.1} we find
\begin{equation}
\begin{aligned}\label{5.2} ie^{i\theta_\gamma}\theta'_\gamma&=(\gamma_1
\gamma'_2-\gamma'_1{\gamma}_2)'
\\& =\gamma_1 (ik_\gamma\gamma'_2-{\gamma}_2)-(ik_{\gamma}
\gamma'_1-\gamma_1)\gamma_2
\\&=i k_\gamma (\gamma_1 \gamma'_2-\gamma'_1{\gamma}_2)
=i k_\gamma e^{i\theta_\gamma},\end{aligned}\end{equation}
which implies statement (2).
\end{proof}

Similarly, we define the {\em Legendre angle}  $\theta_\alpha$ of a
unit speed Legendre curve
$\alpha$ in $H^3_1$ by
\begin{align} \label{5.3}
e^{i\theta_\alpha}=\det{}_{\bf C}(\alpha,\alpha')=\alpha_1
\alpha'_2-\alpha'_1{\alpha}_2 . %\;\; ({\rm mod}\;2\pi).
\end{align}
For instance, the Legendre angle of $\alpha(t)=(\sinh t,\cosh t)$ is 0 (mod $2\pi$).

We also have the following.
\begin{lemma} \label{L:2}  Let $\alpha:I_2\to
H^3_1\subset {\bf C}^2$ be a unit speed curve. Then we have
\begin{enumerate}
\item[(1a)] $\alpha$ is a solution of the second order differential equation$:$
\begin{align}\label{3.8} \alpha''-ik_\alpha \alpha'-\alpha=0,
\end{align}
where $k_\alpha$ is the curvature of $\alpha$ in $H^3_1$.
\item[(1b)] If $\alpha $ satisfies \eqref{3.8}, then $\alpha $ is a Legendre curve if and only if 
$\left<\alpha' (0),i\alpha (0)\right>=0$ $(0\in I_2)$.
\item[(2)] The Legendre angle $\theta_\alpha$ satisfies ${\theta}'_\alpha=k_\alpha$.
\end{enumerate}
\end{lemma}

\begin{proof} This can be done in the same way as Lemma \ref{L:1}. \end{proof}

\begin{remark}
If $(\gamma_1,\gamma_2)$ is a Legendre curve in $S^3$, $(e^{i\theta}\gamma_1,\gamma_2)$ and $(\gamma_1,e^{i\theta}\gamma_2)$ are also Legendre curves in $S^3$. The same happens to a Legendre curve $(\alpha_1,\alpha_2)$ in $H_1^3$.

Using this fact, up to congruences in ${\bf C}^2$,
we can restrict our attention in our construction \eqref{2.3} of Lagrangian surfaces to consider the  initial conditions 
\begin{align}\label{icgamma}
\gamma (0)=(\cos \psi, \sin \psi), \, \gamma ' (0)= e^{i a} (\sin \psi, -\cos \psi), \, 
0\leq \psi \leq \pi/ 2, \, -\pi < a \leq \pi ,
\end{align}
and 
\begin{align}\label{icalpha}
\alpha (0)= (\sinh \delta, \cosh \delta), \, \alpha'(0)=e^{i b}(\cosh \delta, \sinh \delta), \,
\delta \geq 0, \, -\pi < b \leq \pi.\end{align}
\end{remark}

\vspace{.2cm}

We consider now the  Lagrangian conformal immersion $\phi:I_1\times I_2\to {\bf C}^2$ defined by
$\phi(t,s)=\(\alpha_1(t)\gamma_1(s),\alpha_2(t)\gamma_2(s)\)$,
where $\gamma=(\gamma_1,\gamma_2)$ is a unit
speed Legendre curve in  $S^3\subset
{\bf C}^2$  and $\alpha=(\alpha_1,\alpha_2)$  is a unit speed
Legendre curve in $H^3_1\subset {\bf C}^2$. 

The Lagrangian angle map $\beta_\phi$ of $\phi$ (see Section 2) can be computed  by
\begin{align}\label{5.4}e^{i\beta_\phi}={\det}_{\bf
C}(\phi_* e_1,\phi_* e_2), %\;\; ({\rm mod}\;2\pi),
\end{align}
where $e_1,e_2$ is any oriented orthonormal basis of the Lagrangian surface.

We now prove a useful interesting
additive formula which relates the Lagrangian angle map of our
Lagrangian surfaces with the Legendre angle of the generating curves.

\begin{theorem} \label{T:2} Let $\gamma$ be a unit speed
Legendre curve in $S^3$ and $\alpha$ a unit speed Legendre curve in $H^3_1$.
Then  the Lagrangian angle map $\beta_\phi$ of the  Lagrangian conformal immersion
$\phi(t,s)=\(\alpha_1(t)\gamma_1(s),\alpha_2(t)\gamma_2(s)\)$
    and the
Legendre angles $\theta_\gamma$ and $\theta_\alpha$ of
$\gamma$ and $\alpha$ are related by
\begin{align}\label{5.6}\beta_{\phi}(t,s)=
\theta_{\gamma}(s)+\theta_{\alpha}(t)+\pi\;\; \;({\rm mod}\;
2\pi).\end{align}
\end{theorem}
\begin{proof} From  \eqref{2.1}, \eqref{2.2}, \eqref{5.1} and
\eqref{5.3}, we have
\begin{equation}
\begin{aligned}\label{5.7}
e^{i(\theta_\gamma+\theta_\alpha)}&=(\gamma_1
\gamma'_2-\gamma'_1 \gamma_2)(\alpha_1
\alpha'_2-\alpha'_1 \alpha_2)
\\&= \gamma_1\gamma'_2\alpha_1\alpha'_2-\gamma'_1\gamma_2\alpha_1\alpha'_2
-\gamma_1\gamma'_2\alpha'_1\alpha_2
+\gamma'_1\gamma_2\alpha'_1 \alpha_2
\\&= \gamma_1\gamma'_2|\alpha_1|^2\frac{\alpha_1'}{\bar\alpha_2}
+\frac{\gamma'_2}{\bar\gamma_1} |{\gamma}_2|^2|\alpha_1|^2
\frac{\alpha_1'}{\bar\alpha_2}-\gamma_1\gamma'_2\alpha'_1\alpha_2
-\frac{\gamma'_2} {\bar\gamma_1} |\gamma_2|^2\alpha'_1 \alpha_2
\\&=\frac{\alpha_1'\gamma'_2}{\bar\alpha_2\bar\gamma_1}
\(|\gamma_1|^2|\alpha_1|^2 +|\gamma_2|^2|\alpha_1|^2
-|\gamma_1|^2|\alpha_2 |^2 -|\gamma_2|^2|{\alpha}_2|^2\)
\\&=-\frac{\alpha_1'\gamma'_2}{\bar\alpha_2\bar\gamma_1}.\end{aligned}\end{equation}

On the other hand, from \eqref{2.1}, \eqref{2.2}
and \eqref{5.4} we find
\begin{equation}\begin{aligned}\label{5.8}
(|\alpha_1|^2+|\gamma_1|^2)e^{i\beta_\phi}&= {\alpha_1'\alpha_2\gamma_1\gamma'_2
-\alpha_1\alpha'_2\gamma'_1\gamma_2}
\\&= {\alpha_1'\alpha_2\gamma_1\gamma'_2+\alpha_1\frac{\alpha_1'\bar\alpha_1}{\bar\alpha_2}
\frac{\gamma'_2\bar{\gamma}_2}{\bar{\gamma}_1}{\gamma}_2}
\\&=\frac{\alpha_1'\gamma'_2}{\bar\alpha_2\bar\gamma_1}\(|\alpha_2|^2|\gamma_1|^2 +|\alpha_1|^2|\gamma_2|^2\)\\&=\frac{\alpha_1'\gamma'_2}{\bar\alpha_2\bar\gamma_1}
\(\(1+|\alpha_1|^2\)|\gamma_1|^2+|\alpha_1|^2\(1-|\gamma_1|^2\)\)
\\&=\frac{\alpha_1'\gamma'_2}{\bar\alpha_2\bar\gamma_1}\(|\alpha_1|^2+|\gamma_1|^2\),
\end{aligned}\end{equation}
which implies that
\begin{equation}\begin{aligned}\label{5.9} e^{i\beta_\phi}
=\frac{\alpha_1'\gamma'_2}{\bar\alpha_2\bar\gamma_1}.
\end{aligned}\end{equation}

Combining \eqref{5.7} and \eqref{5.9} yields \eqref{5.6}.
\end{proof}

\begin{corollary} \label{C:1} Let $\gamma$ be a unit speed Legendre curve in $S^3$ and $\alpha$ be a unit speed Legendre curve in $H^3_1$. Consider the  Lagrangian conformal
immersion $\phi:I_1\times I_2\to {\bf C}^2$ defined by
$\phi(t,s)=\(\alpha_1(t)\gamma_1(s),\alpha_2(t)\gamma_2(s)\)$.
Then the mean curvature vector field of $\phi$ is given by
\begin{align}\label{4.2} H=\frac{e^{-2u}}{2}\(k_\alpha
J\phi_t+k_\gamma J\phi_s\),\end{align}
where $e^{2u}=|\gamma_1|^2+|\alpha_1|^2$ and $k_\alpha $ and
$k_\gamma$ are the curvature functions of $\alpha $ and $\gamma$,
respectively.
\end{corollary}

\begin{proof} 
According to \eqref{beta}, we have to compute the gradient of the Lagrangian angle $\beta_\phi$. If $e_1:=e^{-u}\partial_t$ and $e_2:=e^{-u}\partial_s$, it is clear that $e_1(\beta_\phi)=e^{-u}k_\alpha$ and $e_2(\beta_\phi)=e^{-u}k_\gamma$ and so \eqref{4.2} follows immediately.
\end{proof}

\section{Applications}

In this section we are devoted to study several families of Lagrangian surfaces of our construction;  those characterized by different geometric properties related with the behaviour of the mean curvature vector.

\subsection{Minimal Lagrangian immersions}

As the first consequence we can obtain from Corollary \ref{C:1} is the following.

\begin{theorem} \label{T:1} Let $\gamma$ be a unit speed
Legendre curve in $S^3$ and let $\alpha$ be a unit speed
Legendre curve in $H^3_1$.
Then the  Lagrangian conformal
immersion $\phi$ defined by
$\phi(t,s)=\(\alpha_1(t)\gamma_1(s),\alpha_2(t)\gamma_2(s)\)$
is minimal if and only if the
Legendre curves
$\gamma$ and $\alpha$ are geodesics in $S^3$ and $H^3_1$,
respectively.
\end{theorem}

Theorem \ref{T:1} also follows directly from Theorem \ref{T:2} by using Lemmas
\ref{L:1}(2) and \ref{L:2}(2) and by taking into account that the
minimality of $\phi $ is equivalently to the constancy of $\beta_\phi$.

Using the statements (1) of Lemmas \ref{L:1} and \ref{L:2}, the
unit speed Legendre curves that are 
geodesic of $S^3$ and $H^3_1$ can be written as $\gamma(s)=\cos s \,
\gamma(0) + \sin s \, \gamma'(0)$ and $\alpha(t)=\cosh t \,\alpha (0)
+ \sinh t \, \alpha '(0)$. After choosing the initial
conditions given in \eqref{icgamma} and \eqref{icalpha}, we arrive at the explicit expressions of the minimal Lagrangian
surfaces in ${\bf C}^2$ that can be constructed by our method taking
\begin{align}\label{geodesic1}
\gamma(s)=(c_\psi \cos s + e^{i a} s_\psi \sin s,s_\psi \cos s - e^{i a} c_\psi \sin s),
\end{align}
where $c_\psi :=\cos \psi$ and $s_\psi := \sin \psi$,
and
\begin{align}\label{geodesic2}
\alpha(t)=(sh_\delta \cosh t + e^{i b} ch_\delta \sinh t,ch_\delta \cosh t + e^{i b} sh_\delta \sinh t),
\end{align}
where $sh_\delta :=\sinh \delta$ and $ch_\delta :=\cosh \delta$.

The Legendre geodesics \eqref{geodesic1} project by the Hopf fibration in the great circles of $S^2(1/2)$ contained in the planes $c_a x_2 = s_a (s_{2\psi}x_3 -c_{2\psi} x_1)$. 
The Legendre geodesics \eqref{geodesic2} project by the Hopf fibration in the geodesics of $H^2(-1/2)$ contained in the planes $c_b x_2 = s_b (ch_{2\delta}x_1 -sh_{2\delta} x_3)$.
So, if $(a,\psi)\in \{ (0,0),(0,\pi/2),(\pi,0),(\pi,\pi/2) \}$ and $(b,\delta) \in \{ (0,0), (\pi,0) \}$, we arrive at
totally geodesic Lagrangian planes.

Topologically all these surfaces are ${\bf R} \times S^1$
and it is possible to prove that they correspond, after changing 
suitably  the complex structure in ${\bf C}^2$ (cf \cite{CM}), to the 
family of  complex surfaces in ${\bf C}^2$ with finite total 
curvature $-4\pi$  (including the Lagrangian catenoid of \cite{cu2}) 
given by Hoffman and Osserman in Proposition 6.6, case 2, of \cite{HO}.

In conclusion, if  we choose great circles in $S^2(1/2)$ and
geodesics  in $H^2(-1/2)$, our construction provides us examples of
minimal Lagrangian surfaces in ${\bf C}^2$.

\subsection{Lagrangian surfaces with constant mean curvature}
 
The easiest (non minimal) examples of Lagrangian surfaces with constant mean curvature, i.e., $|H|\equiv \rho >0$, are those with parallel mean curvature vector.

\begin{theorem} \label{T:new1} Let $\gamma$ be a unit speed
Legendre curve in $S^3$ and let $\alpha$ be a unit speed
Legendre curve in $H^3_1$.
Then the  Lagrangian conformal
immersion $\phi$ defined by
$\phi(t,s)=\(\alpha_1(t)\gamma_1(s),\alpha_2(t)\gamma_2(s)\)$
has parallel $($non null$\,)$ mean curvature vector if and only if the
Legendre curves
$\gamma$ and $\alpha$ in $S^3$ and $H^3_1$ respectively, satisfy that $|\gamma_1|$ and $|\alpha_1|$ are constant.
\end{theorem}

\begin{proof} Using Corollary \ref{C:1}, it is easy to check that  $JH=-(e^{-2u}/{2})(k_\alpha \partial_t+k_\gamma \partial_s )$ is a parallel vector field if and only if  
\begin{align}\label{parallelH} k_\alpha' -u_t k_\alpha +u_s k_\gamma =0, \, 
u_t k_\gamma +u_s k_\alpha =0, \,k_\gamma'-u_s k_\gamma +u_t k_\alpha =0,
\end{align}
where $e^{2u}=|\gamma_1|^2+|\alpha_1|^2$.
From the first and third equation of \eqref{parallelH} we deduce that $k_\alpha'+k_\gamma'=0$ and so 
$k_\alpha (t)=at+b$ and $k_\gamma (s)= -as+c$, with $a,b,c\in {\bf R}$.
We distinguish three cases: We first suppose that $u_s=0$. It is equivalent to $|\gamma_1|$ is constant. Using  \eqref{parallelH} and that $\phi $ is non minimal, we obtain that $u_t=0$ what means that $|\alpha_1|$ is constant.
If $u_t=0$, we make a similar reasoning. Finally, if $u_t\neq 0$ and $u_s \neq 0$, from the second equation of  \eqref{parallelH}, there exists $c_1\in {\bf R^*}$ such that $k_\gamma = -c_1 u_s$ and $k_\alpha = c_1 u_t$. So $u_s=(as-c)/c_1$ and $u_t=(at+b)/c_1$. Putting this in \eqref{parallelH}, we arrive at $a=b=c=0$, which is a contradiction.
\end{proof}

If we call a small circle $\xi$ in $S^2(1/2)$ (respectively, in $H^2(-1/2)$) {\it horizontal} when it is orthogonal
to the $x_3$-coordinate, then we can  easily show that a unit speed
Legendre curve $\gamma$ in $S^3$ is a horizontal lift of a horizontal
circle in $S^2(1/2)$ if and only if $|\gamma_1|$  is a
nonzero constant. Moreover, such Legendre curves can be parametrized by
\begin{equation}\label{2.21}\gamma(s)= \left(\cos \psi \, e^{i \tan \psi \, s}, \sin \psi  \, e^{-i \cot \psi \, 
s}\right), \  \psi \in (0,\pi/2),\end{equation}
where $\pi/2 - 2\psi $ is the latitude of the parallel $\pi \circ \gamma $.

Similarly, a unit speed Legendre curve $\alpha$ in $H_1^3$ is a horizontal lift
of a horizontal circle in $H^2(-1/2)$ if and only if $|\alpha_1|$  is
a nonzero constant. Moreover, such Legendre curves can be parametrized by
\begin{equation}\label{2.30}\alpha(t)=\left(  \sinh \delta \, e^{i\coth \delta \, t}, \cosh \delta  \, e^{i  \tanh \delta \, t }\right),\  \delta >0.\end{equation}

In conclusion, using both Legendre curves given in \eqref{2.21} and in \eqref{2.30} 
in our construction we obtain conformal parametrizations of the examples of Lagrangian 
surfaces with parallel mean curvature vector in ${\bf C}^2$. They correspond to flat tori $S^1 \times S^1$ in the 3-sphere of radius $\sqrt{\sin ^2 \psi + \sinh ^2 \delta}$.

\begin{theorem} \label{T:new2} Let $\gamma$ be a unit speed Legendre curve in $S^3$ and let $\alpha$ be a unit speed Legendre curve in $H^3_1$. Then the  Lagrangian conformal
immersion $\phi$ defined by
$\phi(t,s)=\(\alpha_1(t)\gamma_1(s),\alpha_2(t)\gamma_2(s)\)$
has constant mean curvature $|H|\equiv \rho >0$ if and only if the
Legendre curves
$\gamma$ and $\alpha$ in $S^3$ and $H^3_1$, respectively, satisfy that $k_\gamma^2=4\rho^2 |\gamma_1|^2- \lambda$ and $k_\alpha^2 = 4\rho^2 |\alpha_1|^2 + \lambda$ with $\lambda \in {\bf R}$.
\end{theorem}

\begin{proof}
Using Corollary \ref{C:1}, we have that 
$4\rho ^2 (|\gamma_1|^2+|\alpha_1|^2) = k_\alpha ^2+k_\gamma ^2$. Since $\gamma$ depends on $s$ and $\alpha $ depends on $t$, we obtain the result.
\end{proof}

Now we see how the condition on $\gamma $ and $\alpha $ in Theorem \ref{T:new2} determine both curves. Let define
$r(s):=|\gamma_1(s)|$.  Using the Legendre character of $\gamma $ and that it is parametrized by arclength and satisfies
\eqref{3.2}, it is not difficult to get that $\gamma$ can be expressed in terms of $r$ in the following way: 
\begin{align}\label{expgamma}
\gamma (s) =\Bigg (r(s)\exp \text{\small $\Bigg(i \int_0^s \frac{\sqrt{1-r^2-r'^2}}{r}ds\Bigg)$},
              \sqrt{1-r(s)^2}\exp \text{\small $\Bigg({i\int_0^s \frac{r\sqrt{1-r^2-r'^2}}{r^2-1}}ds\Bigg)$}\Bigg).
\end{align}
We observe that when $r(s)=\cos s$, we get the geodesic $\gamma (s) = (\cos s , \sin s)$, and if we take $r$ constant, say
$r\equiv \cos \psi$, we arrive at the expression of \eqref{2.21}. We can also compute the curvature of $\gamma $ in terms of $r=r(s)$ obtaining
\begin{align}\label{k1}
r''-\frac{1-r^2-r'^2}{r}+r+ k_\gamma \sqrt{1-r^2-r'^2}=0.
\end{align}

If we use a similar argument, a unit speed Legendre curve in $H^3_1$ can be written in terms of $r(t)=|\alpha_1(t)|$ as
\begin{align}\label{expalpha}
\alpha (t) = \Bigg(r(t)\exp \text{\small $\Bigg(i \int_0^t \frac{\sqrt{1+r^2-r'^2}}{r}dt\Bigg)$},
              \sqrt{1+r(t)^2}\exp \text{\small $\Bigg({i\int_0^t \frac{r\sqrt{1+r^2-r'^2}}{1+r^2}}dt\Bigg)$}\Bigg).
\end{align}
We note that $r(t)=\sinh t$ give us the geodesic $\alpha (t)=(\sinh t, \cosh t)$ and that $r(t)\equiv \sinh \delta$ leads to \eqref{2.30}. Moreover, the curvature of $\alpha $ is given by
\begin{align}\label{k2}
r''-\frac{1+r^2-r'^2}{r}-r+ k_\alpha \sqrt{1+r^2-r'^2}=0.
\end{align}

We study the case of Theorem \ref{T:new2}. If $k_\gamma(s)^2=4\rho^2 r(s)^2-\lambda$ and $k_\alpha(t)^2=4\rho^2 r(t)^2+\lambda$, we are able to obtain first integrals of the differential equations \eqref{k1} and \eqref{k2}:
\begin{align}\label{4.7}
\frac{(4\rho^2 r^2-\lambda)^{3/2}}{12 \rho^2}+ \mu_1 = r\sqrt{1-r^2-r'^2}, \, r=r(s),
\end{align}
and
\begin{align}\label{4.9} 
\frac{(4\rho^2 r^2+\lambda)^{3/2}}{12 \rho^2}+ \mu_2 = r\sqrt{1+r^2-r'^2}, \, r=r(t),
\end{align} where $\lambda,\mu_1,\mu_2$ are arbitrary constants.

This shows that the family of Lagrangian surfaces with constant mean curvature $\rho >0$ in our construction with Legendre curves is quite big. In general, the solutions of \eqref{4.7} and \eqref{4.9} are not easy to control, appearing hyperelliptic functions in most cases. We finish this section considering the following illustrative situation. 

\medskip

{\bf Particular case}: Let $\lambda =\mu_1 = \mu_2=0$. Up to dilations, we can suppose $\rho = 3/2$. Then equations \eqref{4.7} and \eqref{4.9} reduce to $r'^2+r^2+r^4=1$ and $r'^2-r^2+r^4=1$,  respectively.
After solving the  differential equation $r'^2+r^2+r^4=1$, we know that, up to  translations on $s$, its solution   is given by
\begin{align}\label{4.10}
r(s)=\sqrt{\tfrac{\sqrt 5 -1}{2}} \,\cn\left(\sqrt[\uproot{2} 4]{5}s,k\right), \;\; \; k=\sqrt{\tfrac{5-\sqrt{5}}{10}},
\end{align} 
where $\cn$ is a Jacobi elliptic function usually known as the cosine amplitude and $k$ is its  modulus  (cf., for instance, \cite{L}).

Hence,  using standard formulae on elliptic functions and a straightforward long computation, \eqref{expgamma} and \eqref{4.10} imply that, up to rotations, $\gamma$ is given by 
\begin{equation}\begin{aligned}\label{expgamma.1}
&\gamma (s) =\left( \dn\left(\sqrt[\uproot{2} 4]{5}s,k\right) + i \,  k\, \sn\left(\sqrt[\uproot{2} 4]{5}s,k\right) \right)
\Bigg( \sqrt{\tfrac{\sqrt 5 -1}{2}} \cn\left(\sqrt[\uproot{2} 4]{5}s,k\right),
\\& \sqrt{1+\tfrac{1}{2}(1-\sqrt{5})\cn^2\big(\sqrt[\uproot{2} 4]{5}s,k\big)}\,
\frac{2 \dn\left(\sqrt[\uproot{2} 4]{5}s,k\right) - \sqrt{5+2\sqrt{5}}
 i \,  \sn\left(\sqrt[\uproot{2} 4]{5}s,k\right)}{\sqrt{4 \dn^2\left(\sqrt[\uproot{2} 4]{5}s,k\right) + (5+2\sqrt 5) \sn^2\left(\sqrt[\uproot{2} 4]{5}s,k\right)}}
\Bigg),
\end{aligned}\end{equation}
where $\dn$ and $\sn$ are the Jacobi elliptic function known as the delta amplitude and the sine amplitude with modulus $k$.

Similarly, up to translations in $t$, the solution of $r'^2-r^2+r^4=1$ is given by
\begin{align}\label{4.12}
r(t)=\sqrt{\tfrac{\sqrt 5 +1}{2}} \,\cn\left(\sqrt[\uproot{2} 4]{5}t,\hat k\right),\;\; \hat k=\sqrt{\tfrac{1}{10}(5+\sqrt{5}\,)}.\end{align}

Thus, in an analogous way, it follows from \eqref{expalpha}, \eqref{4.12} and a long computation that, up to rotations, $\alpha$ is given by 
\begin{equation}\begin{aligned}\label{expalpha.1}
&\alpha (t) = \left( \dn\left(\sqrt[\uproot{2} 4]{5}t,\hat k\right) + i \,  \hat k\, \sn\left(\sqrt[\uproot{2} 4]{5}t,\hat k\right) \right)
\Bigg( \sqrt{\tfrac{\sqrt 5 +1}{2}} \cn\left(\sqrt[\uproot{2} 4]{5}t,\hat k\right),
\\& \sqrt{1+\tfrac{1}{2}(1+\sqrt{5})\cn^2\big(\sqrt[\uproot{2} 4]{5}t ,\hat k\big)}\,
\frac{2 \dn\left(\sqrt[\uproot{2} 4]{5}t,\hat k\right) - \sqrt{5-2\sqrt{5}}
 i \,  \sn\left(\sqrt[\uproot{2} 4]{5}t,\hat k\right)}{\sqrt{4 \dn^2\left(\sqrt[\uproot{2} 4]{5}t,\hat k\right) + (5-2\sqrt 5) \sn^2\left(\sqrt[\uproot{2} 4]{5}t,\hat k\right)}}
\Bigg).
\end{aligned}\end{equation}

We remark that both Legendre curves $\gamma $ and $\alpha $ given in \eqref{expgamma.1} and \eqref{expalpha.1} are periodic on account of the periodicity of the elliptic functions $\cn $, $\sn $ and $\dn $.
So they provide an interesting example of a Lagrangian torus with constant mean curvature in complex Euclidean plane.

\subsection{New examples of Hamiltonian-minimal Lagrangian surfaces}

By applying Theorem \ref{T:2}, we have the following

\begin{theorem} \label{T:3} Let $\gamma$ and $\alpha$ be  unit speed
Legendre curves in $S^3$ and  $H^3_1$, respectively.
Then the  Lagrangian conformal immersion
$\phi(t,s)=\(\alpha_1(t)\gamma_1(s),\alpha_2(t)\gamma_2(s)\)$
     is Hamiltonian-minimal if and only if the
curvature functions $k_\alpha$ and $k_\gamma$ of  $\alpha$ and
$\gamma$ are given by $k_\alpha (t)=at+b$ and $k_\gamma (s)= -as+c$ with $a,b,c\in {\bf R}$.
\end{theorem}

\begin{proof} It is known that the Lagrangian surface $\phi$
is Hamiltonian-minimal if and only if the Lagrangian angle map
$\beta_\phi$ is harmonic, i.e., $\Delta\beta_\phi=0$ (see Section 2).
Thus,  Theorem \ref{T:2} implies that $\phi$
is Hamiltonian-minimal if and only if we have
$\theta_\gamma'' + \theta_\alpha''=0$. Therefore, by Lemmas
\ref{L:1} and \ref{L:2}, we know that
$\phi$ is Hamiltonian-minimal if and only if the
curvature function $k_\alpha$ and
$k_\gamma$ of  $\alpha$ and $\gamma$ satisfy $k_\alpha' +k_\gamma ' =0$.
We put then $k_\alpha'=a=-k_\gamma$ and the proof is finished.
\end{proof}

We distingish two essential cases in this family:
\smallskip

{\bf Case (i)}: $a=0$, i.e., $k_\alpha $ and $k_\gamma$ are constant.

  From Lemma \ref{L:1},(1) we know that unit speed Legendre curves in $S^3$ with
constant curvature $k_\gamma \equiv c$ can be parametrized by
$$
\gamma(s)=e^{i(c+\sqrt{c^2+4})s/2}A_1+e^{i(c-\sqrt{c^2+4})s/2}B_1
$$
for suitable $A_1,B_1\in {\bf C}^2$ that can be expressed in terms on the initial conditions given in \eqref{icgamma}.

Similarly,  from Lemma \ref{L:2},(1) we also know that unit speed Legendre 
curves in $H^3_1$ with constant curvature $k_\alpha\equiv b$ can be
parametrized by

\begin{enumerate}
\item If $|b|>2$, $$\alpha(t)=e^{i\, b t/2}  ( e^{i\sqrt{b^2-4}\,t/2}A_2+e^{-i\sqrt{b^2-4}\,t/2}B_2);$$
\item if $|b|<2$, $$\alpha(t)=e^{i\, b t/2}  (e^{\sqrt{4-b^2}\,t/2}A_2+e^{-\sqrt{4-b^2}\,t/2}B_2);$$
\item if $b=2$, $$\alpha(t)=e^t A_2 + t e^t B_2;$$
\item if $b=-2$, $$\alpha(t)=e^{-t} A_2 + t e^{-t} B_2,$$
\end{enumerate}
for suitable $A_2,B_2\in {\bf C}^2$ that can be expressed in terms on the initial conditions given in \eqref{icalpha}.

\begin{remark}
In this context, it is not difficult to check that the Hamiltonian-minimal Lagrangian tori of \cite{cu1} are constructed by using horizontal small circles in $S^2(1/2)$ (see \eqref{2.21}) with
closed horizontal lift (i.e., $\tan^2 \psi $  is a rational number) and the projections to $H^2(-1/2)$ of the
above $\alpha$'s for certain $b$'s such that $|b|>2$.
\end{remark}

\smallskip
{\bf Case (ii)}: $a\neq 0$, i.e., $k_\alpha $ and $k_\gamma$ are certain linear functions of the arc parameter.

In this case, after applying suitable translations,  we have $\kappa_\alpha=at$ and $\kappa_\gamma=-as$. Thus, by Lemmas \ref{L:1} and \ref{L:2}, we know that the Legendre curves $\alpha$ and $\gamma$ satisfy
\begin{align}  \alpha''(t)-i at\alpha'(t)-\alpha(t)=0,\;\; \gamma''(s)+ias\gamma'(s)+\gamma(s)=0.\end{align}
Therefore, after solving these differential equations, we  know that the unit speed Legendre 
curves $\alpha$ in $H^3_1$ with  curvature $\kappa_\alpha=at$ and the unit speed Legendre curves $\gamma$ in $S^3$ with
 $\kappa_\gamma=-as$  can be expressed in terms of Hermite polynomials and hypergeometric functions (see \cite{S}) by
\begin{align} &\alpha(t)= \hbox{\rm HermiteH}\big(i/a,\sqrt{{ia}/{2}}\,t\big)A_1+{}_1F_1\big(1/({2ai}),1/2,ait^2/2\big)B_1,\\&
\gamma(s)=e^{- ais^2/2}\big\{ \hbox{\rm HermiteH}\big(1/(ai)-1,\sqrt{{ia}/2}\,s \big)A_2\\&\hskip.6in \notag + {}_1F_1\big((ai-1)/(2ai),1/2,ais^2/2\big)B_2\big\}
\end{align} for suitable $A_1,B_1\in {\bf C}^2_1$ and $A_2,B_2\in {\bf C}^2$ depending on the initial conditions given in \eqref{icgamma} and \eqref{icalpha}, where HermiteH is the Hermite polynomial  and ${}_1F_1$ is the Kummer confluent hypergeometric function.

\subsection{Willmore Lagrangian surfaces}

Consider the Willmore functional
\begin{align}\label{6.1} W=\int_\Sigma |H|^2dA\end{align}
for a surface $\Sigma $  in a  Euclidean space.

For a unit speed Legendre curve $\gamma$ in $S^3$ and a unit speed Legendre curve  $\alpha$ in $H^3_1$, the Willmore functional of the  Lagrangian conformal immersion $\phi:I_1\times I_2\to {\bf C}^2; (t,s)\mapsto (\alpha_1(t)\gamma_1(s),\alpha_2(t)\gamma_2(s))$
is given by
\begin{align}\label{6.2} W_\phi=\frac{1}{4}\int_{\phi(I_1\times I_2)} |\nabla\beta_\phi|^2dA.
\end{align}
Hence, it follows from Theorem \ref{T:2} and Lemmas \ref{L:1} and \ref{L:2}  that the
Willmore functional associated with $\phi$ is given by
\begin{equation}\begin{aligned}\label{6.3} W_\phi&=\frac{1}{4}\int_{I_1\times I_2}
\(k_\alpha^2+k_\gamma^2\)dtds\\& =\frac{L (\gamma)}{4}\int_{I_1} k_\alpha^2dt +\frac{L(\alpha)}{4}\int_{I_2}
k_\gamma^2ds, \end{aligned}\end{equation}
   where $L(\gamma)$ and $L(\alpha)$
denote the length of $\gamma$ and of $\alpha$, respectively.

\begin{theorem} \label{T:4}  Let $\gamma$ and $\alpha$ be  unit speed Legendre curves in $S^3$ and  $H^3_1$, respectively. Then the Lagrangian conformal immersion
$\phi(t,s)=\(\alpha_1(t)\gamma_1(s),\alpha_2(t)\gamma_2(s)\)$    is a  critical point of the Willmore functional $W_\phi$ {\rm (}with fixed lengths $L(\alpha) $ and $L(\gamma))$  if  and only if the Legendre curves
$\alpha$ and $\gamma$ are elastic curves. \end{theorem}

\begin{proof} From \eqref{6.3}, we see that the critical points
of the Willmore functional $W_\phi $ (with fixed $L_1=L(\alpha) $ and
$L_2=L(\gamma)$) are given by the Lagrangian conformal immersions  constructed with
Legendre curves $\alpha$ and $\gamma$ that are critical points of the
functionals $\int_0^{L_1} k_\alpha^2dt$ and $\int_0^{L_2} k_\gamma^2ds$, respectively.
But these are precisely free elastic curves according to \cite{LS}. \end{proof}

\begin{remark}
As corollary of Theorem \ref{T:4}, using free elastica in $S^2(1/2)$ 
and $H^2(-1/2)$ our construction provides new examples of Willmore 
Lagrangian surfaces in ${\bf C}^2$.  In \cite{LS} we can find 
explicitly examples of (closed) free elastica on the sphere and in 
the Poincar\'e disk.

A different construction of Willmore Lagrangian surfaces in ${\bf C}^2$ can be found in \cite{cu3}.
\end{remark}

\vskip.2in 
\noindent {\small Departamento de Matem\'aticas\\ Escuela
Polit\'ecnica Superior Universidad de Ja\'en\\  23071, Ja\'en\\
Spain}

\noindent {\small {\it E-mail address}: {icastro@ujaen.es} }

\vskip.1in 
\noindent {\small
{Department of Mathematics\\
    Michigan State University\\ East Lansing\\ Michigan 48824--1027\\ U.S.A.}
    
\noindent {\small {\it E-mail address}:  {bychen@math.msu.edu}}

\end{document}